\documentclass[journal]{IEEEtran}
\usepackage{cite}

\usepackage{color}

\usepackage[pdftex]{graphicx}

\usepackage{amsmath}
\usepackage{amsfonts}
\usepackage{amsthm}

\theoremstyle{plain}

\newtheorem{remark}{Remark}

\usepackage{algorithm}
\usepackage{algorithmic}
\usepackage{multirow}  
\usepackage{booktabs}  

\usepackage{array}

\ifCLASSOPTIONcompsoc
  \usepackage[caption=false,font=normalsize,labelfont=sf,textfont=sf]{subfig}
\else
  \usepackage[caption=false,font=footnotesize]{subfig}
\fi
\usepackage{dblfloatfix}

\renewcommand{\arraystretch}{1.3}

\usepackage{nomencl}
\makenomenclature

\begin{document}
%
\title{A Chance Constrained Information-Gap Decision Model for Multi-Period Microgrid Planning}
%
%
%

\author{Xiaoyu Cao,~\IEEEmembership{Student Member,~IEEE,}
        Jianxue Wang,~\IEEEmembership{Member,~IEEE,}
        and~Bo~Zeng,~\IEEEmembership{Member,~IEEE}
\thanks{Xiaoyu Cao and Jianxue Wang are with the School
of Electrical Engineering, Xi'an Jiaotong University, Xi'an,
Shaanxi 710049 China (e-mail: cxykeven2013@stu.xjtu.edu.cn; jxwang@mail.xjtu.edu.cn.)}
\thanks{Bo Zeng is with the Department of Industrial Engineering and the Department of Electrical and Computer Engineering, University of Pittsburgh, Pittsburgh,
PA 15106 USA (e-mail: bzeng@pitt.edu).}}

\maketitle

\begin{abstract}
This paper presents a chance constrained information gap decision model for multi-period microgrid
expansion planning (MMEP) considering two categories of uncertainties, namely random and non-random
uncertainties. The main task of MMEP is to determine the optimal sizing, type selection, and installation time of distributed
energy resources (DER) in microgrid. In the proposed formulation, information gap decision theory (IGDT) is applied to hedge against
 non-random uncertainties of long-term demand growth. Then, chance constraints are imposed in the operational stage
to address the random uncertainties of hourly renewable energy generation and load variation. The objective of chance constrained information gap decision model is to maximize the robustness level of DER investment meanwhile satisfying a set of operational constraints with a high probability. The integration of IGDT and chance constrained program, however, makes it very challenging to compute. To address this challenge, we propose and implement a strengthened bilinear Benders decomposition method. Finally, the effectiveness of proposed planning model is verified through the numerical studies on both the simple and practical complex microgrid. Also, our new computational method demonstrates a superior solution capacity and scalability. Compared to directly using a professional mixed integer programming solver, it could reduce the computational time by orders of magnitude.
\end{abstract}


%
\IEEEpeerreviewmaketitle

\nomenclature{$\alpha_{\rm L}$}{Horizon of load uncertainty}%
\nomenclature{$k$}{Index of DER}%
\nomenclature{$h$}{Index of hour}%
\nomenclature{$d$}{Index of day}%
\nomenclature{$t$}{Index of year}%
\nomenclature{$T$}{Number of years}%
\nomenclature{$D$}{Number of typical days in a year}%
\nomenclature{$H$}{Number of hours in a day}%
\nomenclature{$r$}{Discount rate}%
\nomenclature{$\tau(t)$}{Present value factor in $t$th year}%
\nomenclature{$X_k^t$}{Number of installed DER $k$ in $t$th year}%
\nomenclature{$IC_k$}{Unit cost of DER $k$ (except for battery packs)}%
\nomenclature{$PC_k$}{Unit power cost of ESS $k$}%
\nomenclature{$EC_k$}{Unit energy cost of ESS $k$}%
\nomenclature{$OC_k^t$}{Yearly O\&M rate of DER $k$}%
\nomenclature{$c_k$}{Unit fuel cost of DFG $k$}%
\nomenclature{$q_{\rm v}$}{Penalty cost for load curtailment}%
\nomenclature{$P_k^{\max}$}{Nominal power generation of DER $k$}%
\nomenclature{$E_k^{\max}$}{Nominal storage capacity of ESS $k$}%
\nomenclature{$C_{\rm ann}^{\max}$}{Upper limit of annual investment cost}%
\nomenclature{$X_k^{\max}$}{Upper limit of annual installation of DER $k$}%
\nomenclature{$\psi^t$}{Load capacity of microgrid in $t$th year}%
\nomenclature{$R^t$}{Long-term reserve requirement for microgrid in $t$th year}%
\nomenclature{$RU_k$}{Ramp up rate of DFG $k$}%
\nomenclature{$RD_k$}{Ramp down rate of DFG $k$}%
\nomenclature{$OR^{tdh}$}{Short-term reserve requirement for microgrid at $h$th hour in $d$th day of $t$th year}%
\nomenclature{$P_k^{tdh}$}{Power generation of RES and DG $k$ at $h$th hour in $d$th day of $t$th year}%
\nomenclature{$P_{k,\rm dis}^{tdh}$}{Discharge power of ESS $k$ at $h$th hour in $d$th day of $t$th year}%
\nomenclature{$P_{k,\rm ch}^{tdh}$}{Charge power of ESS $k$ at $h$th hour in $d$th day of $t$th year}%
\nomenclature{$P_{\rm lc}^{tdh}$}{Load curtailment at $h$th hour in $d$th day of $t$th year}%
\nomenclature{$E_k^{tdh}$}{Stored energy of ESS $k$ at $h$th hour in $d$th day of $t$th year}%
\nomenclature{$\eta_k$}{Charge/discharge efficiency of ESS $k$}%
\nomenclature{$p_{\rm L}^{tdh}$}{Load variation factor at $h$th hour in $d$th day of $t$th year}%
\nomenclature{$A_k^{tdh}$}{Available capacity of RES $k$ at $h$th hour in $d$th day of $t$th year}%
\nomenclature{$\omega$}{Minimum installation ratio of DFGs}%
\printnomenclature

\section{Introduction}
%
%
%
%
\IEEEPARstart{M}{icrogrid} has been considered as a viable solution for the integration and management of distributed energy resources (DER)~\cite{Hatziargyriou2007Microgrids}, which typically include renewable energy sources (RES), energy storage system (ESS), dispatchable fuel-generators (DFG), etc. The primary goal of multi-period microgrid expansion planning (MMEP) is to determine the optimal sizing, type selection, and installation time of DERs. The multi-period investment strategy is considered as an effective way to reduce the budget level of microgrid planning~\cite{Hajipour2015Stochastic}, which would benefit the commercial customization of microgrid technology in the long run.

However, it is very challenging to manage the uncertainties in microgrid expansion planning. First, the intermittency of renewable energy generation increases the randomness of nodal power injection~\cite{Khodaei2014Microgrid}. Second, the rapid development of renewable energy industry brings non-trivial uncertainties into the long-term investment of power sources~\cite{Jeon2014Long}. At last, it is difficult to predict the load variation in a microgrid due to the lack of historical data as well as the disaggregated environment of distributed generation system~\cite{Hernandez2014A}. Either underestimating or overrating the impact of uncertainties in planning decisions could possibly jeopardize the power-supply reliability and cost-benefit of a microgrid.

\indent In order to develop an advanced uncertainty management scheme, it is necessary to classify the uncertainties for MMEP problem. Here, we adopt the classical categorization for uncertain factors used in~\cite{Burke1988Trade,Buygi2004Market,Buygi2006Network,Maghouli2011A}:
\begin{itemize}
\item[$\bullet$] \textbf{Random uncertainties} follow the specified probability distributions, which are high-frequency and repeatable in the short period. The hourly variation of load demand and RES generation is a typical example in this category.
\item[$\bullet$] \textbf{Non-random uncertainties} occur in lower frequency and larger time scale, which are difficult to accurately fit into any distribution due to data insufficiency. The uncertainties in demand growth fall into this category.
\end{itemize}

Probabilistic optimization methods, e.g. stochastic program (SP) and chance constrained program (CCP), are popular tools to manage the random uncertainties in microgrid planning. The main idea of SP is to employ a finite set of scenarios to represent the possible realizations of random variables. Ref~\cite{Bahramirad2012Reliability} utilized SP to optimize the one-dispose investment for ESS in microgrid, where the random scenarios of RES generation were involved. Similar method was adopted in~\cite{Hajipour2015Stochastic} for the capacity expansion planning of wind turbines (WT) and ESS in a stand-alone microgrid. Typically, an SP model considers every scenario in the scenario pool, which may lead to a costly solution that is less practical. Unlike SP, the key idea of CCP is to derive a solution that performs well in a subset of scenarios with a probabilistic guarantee, e.g., a solution works well in 90\% random situations. Hence, CCP is a decision-making scheme providing a trade-off between cost and performance. Typical applications of CCP in microgrid planning can be found in~\cite{Guo2014Multi,Liu2011Optimal,Kamjoo2014Chance}. Ref~\cite{Guo2014Multi} presented a simulation based multi-objective CCP model to derive the single-stage planning solution of a stand-alone microgrid. Such simulation-based framework was also adopted in~\cite{Liu2011Optimal} to solve a chance constrained planning problem that addressed the siting and sizing of DERs in a distribution network. Ref~\cite{Kamjoo2014Chance} applied CCP to design a hybrid renewable energy system considering the random generation of WT and photovoltaic arrays (PV). We also note that a couple of other probabilistic methods have also been adopted to handle the randomness in microgrid planning, e.g., Markovian sizing approach~\cite{Dong2016Storage}, sample average approximation method~\cite{Sharafi2015Stochastic}, and conditional value-at-risk (CVaR)~\cite{Mena2014A}, which have different trade-offs between modeling advantages and computational requirements. The common feature of those methods is to make use of the uncertainty information contained in rich data. However, they are less applicable when data is unavailable or inaccurate.

The microgrid planning under non-random uncertainties, whose distributions are unknown, has been investigated using uncertainty set-based approaches, e.g. robust optimization (RO) and information-gap decision theory (IGDT). They only need uncertainty sets rather than probability distributions and make planning decisions under the worst scenarios within those sets. Ref~\cite{Billionnet2016Robust} proposed and computed a two-stage RO based microgrid planning model. Ref~\cite{Khodaei2014Microgrid} studied microgrid planning using a RO scheme to consider the uncertainties embedded in long-term load forecast, volatility of RES generation, and unintentional islanding. {Different from RO, which requires for an uncertainty set with explicit boundaries, IGDT has the expected system performance as an input while the objective is to maximize the horizon of the uncertainty set, in which the expected performance is guaranteed \cite{Ben2006Info}. Hence, for the situation where the system performance target is clear (e.g., the budget plan of microgrid projects), IGDT helps the planners to identify a planning solution with a maximized immunity (subject to the uncertainty budget) to handle future uncertainties. In this regard, IGDT needs even less information about uncertainty structure~\cite{Kazemi2015Risk}, which makes it adaptive for the problems under unstructured uncertainties.} Given its modeling capacity against uncertainties, we observe that IGDT has been applied to the expansion planning of generation system~\cite{Shivaie2016A} and transmission network~\cite{Taherkhani2014IGDT}, but is rarely mentioned in microgrid planning.


Based on the literatures, we note that the choice of planning strategy largely depends on the type of uncertain factors. It becomes much more complicated when the distribution information is partially available, e.g., there is a mixture of random and non-random uncertainties. Such case, which usually happens in practical planning problems, demands for an integrated optimization framework to manage both types of uncertainties. Hence, a chance constrained information-gap decision (IGD) model is proposed in this paper. Our microgrid planning formulation considers multi-period investment scheme and detailed operational modeling to capture the uncertainties under multiple time-scales, i.e., 1) the non-random uncertainty of yearly demand growth; 2) the random uncertainties of hourly RES generation and load variation. We follow the study presented in~\cite{Hajipour2015Stochastic,Guo2014Multi,Kamjoo2014Chance} to derive a planning solution for a stand-alone microgrid. Compared to a grid-supported system, the stand-alone microgrid is typically more vulnerable to uncertain factors, which hence requires us to use advanced optimization methods to incorporate those factors on system design. We note that the proposed model can also be extended to study a grid-connected microgrid by including additional decision variables and operational constraints concerning the power exchange between the microgrid and the utility grid. {{Also, given the close electrical connection between load and DERs in microgrids~\cite{Khodaei2014MicrogridScheduling}, we do not include the network transmission and associated reactive power or power factor issue in our model. Hence, our MMEP formulation can be categorized as a resource planning problem without regard to network representation. }}

Next, we summarize the main contributions of the presented research:

\begin{enumerate}
\item We construct the MMEP model under an integrated optimization framework of IGDT and CCP. In the derived framework, IGDT is applied to maximize the robustness level against non-random uncertainties through multi-period investment strategy, while the chance constraints are imposed in the operational stage to provide a flexible scheme to balance the random uncertainties and operational cost.
\item To compute this integrated model, a strengthened customization of bilinear Benders decomposition method~\cite{Zeng2014Chance} is developed with some non-trivial changes to adapt it for further reducing the computational burden from the complicated combinatorial structure.
\item We implement a series of comparative numerical tests on a simple test microgrid to examine the validity of the proposed planning method given the existence of both random and non-random uncertainties. Also, the scalability of our computational method is verified by a practical complex microgrid.
\end{enumerate}

The remainder of this paper is organized as below. Section 2 provides the basic form of MMEP problem. Section 3 derives the mathematical formulation of chance constrained IGD model for microgrid planning. Section 4 illustrates the solution approach of strengthened bilinear Benders decomposition. Section 5 shows the results of numerical tests. Finally, conclusions are drawn in Section 6.

\section{Multi-Period Microgrid Expansion Planning}
\label{Sec:MMEP}
The basic form of MMEP problem contains a least-cost objective function and a set of constraints from investment and operational considerations, which restrict investment and operational decisions of alternative DERs. Three categories of DERs are considered, i.e., $\rm DER=\{RES, DFG, ESS\}$. Our basic MMEP model, which characterizes the multi-period investment and long-term scheduling schemes, is formulated as:
\setlength{\arraycolsep}{0.1em}
\begin{eqnarray}
&&\min \ TB = \sum_{t=1}^T \tau(t) (C_{\rm inv}^t+C_{\rm opr}^t) \label{OBJ}\\
&&C_{\rm inv}^t=\sum_{k\in \{\rm RES,DFG\}} IC_k P_k^{\max} (X_k^t-X_k^{t-1}) \nonumber\\
&&+ \sum_{k\in{\rm ESS}} (PC_k P_k^{\max}+EC_k E_k^{\max}) (X_k^t-X_k^{t-1}) \quad \forall t \label{OBJ_1}\\
&&C_{\rm opr}^t=\sum_{k\in {\rm DER}} OC_k^t P_k^{\max} (X_k^t-X_k^{t-1}) \nonumber\\
&&+\sum_{d=1}^D \sum_{h=1}^H \frac{365}{D} (\sum_{k\in{\rm DFG}}c_k P_k^{tdh}+q_{\rm v}P_{\rm lc}^{tdh}) \quad \forall t \label{OBJ_2}\\
&&     \nonumber \\
&&s.t. \nonumber \\
&&C_{\rm inv}^t \leq C_{\rm ann}^{\max} \quad \forall t \label{TB}\\
&&X_k^{t-1} \leq X_k^t \quad \forall k, \forall t \label{YearGap}\\
&&0 \leq X_k^t \leq X_k^{\max} \quad \forall k, \forall t  \label{XLimit}\\
&&\sum_{k\in {\rm DER}} P_k^{\max} X_k^t \geq \psi^t + R^t \quad \forall t \label{QB1}\\
&&\sum_{k\in {\rm DFG}} P_k^{\max} X_k^t \geq \omega \psi^t \quad \forall t \label{QB2}\\
&&     \nonumber \\
&&\sum_{k\in \{\rm RES,DFG\}} P_k^{tdh}+\sum_{k\in{\rm ESS}} (P_{k,\rm dis}^{tdh}-P_{k,\rm ch}^{tdh})\geq \psi^t p_{\rm L}^{tdh}-P_{\rm lc}^{tdh} \nonumber\\
&&\qquad \qquad \qquad \qquad \qquad \qquad \quad \quad \quad \quad \quad \quad \forall t, \forall d, \forall h \label{PB}\\
&&0\leq P_k^{tdh} \leq X_k^t A_k^{tdh} \quad \forall k\in{\rm RES}, \forall t, \forall d, \forall h \label{Pres}\\
&&X_k^t P_k^{\min} \leq P_k^{tdh} \leq X_k^t P_k^{\max} \quad \forall k\in{\rm DFG}, \forall t, \forall d, \forall h \label{Pdg}\\
&&-RD_k \leq P_k^{tdh}-P_k^{tdh-1} \leq RU_k \ \forall k\in{\rm DFG}, \forall t, \forall d, \forall h \label{Pramp}\\
&&\sum_{k\in {\rm DFG}} (X_k^t P_k^{\max}-P_k^{tdh}) \geq OR^{tdh} \quad \forall t, \forall d, \forall h \label{Preserve}\\
&&0 \leq P_{k, \rm ch}^{tdh}\leq X_k^t P_k^{\max} \quad \forall k\in{\rm ESS}, \forall t, \forall d, \forall h \label{Pbes1}\\
&&0 \leq P_{k, \rm dis}^{tdh}\leq X_k^t P_k^{\max} \quad \forall k\in{\rm ESS}, \forall t, \forall d, \forall h \label{Pbes2}\\
&&X_k^t E_k^{\min} \leq E_k^{tdh} \leq X_k^t E_k^{\max} \quad \forall k\in{\rm ESS}, \forall t, \forall d, \forall h \label{Ebes}\\
&&E_k^{tdh}=\sum_{i=1}^h (\eta_k P_{k,\rm ch}^{tdi}-P_{k,\rm dis}^{tdi}/\eta_k) \ \forall k\in{\rm ESS}, \forall t, \forall d, \forall h \label{PEbes}\\
&&E_k^{td0}=E_k^{tdH} \quad \forall k\in{\rm ESS}, \forall t, \forall d \label{Einitial}\\
&&X_k^t \in \mathbb{Z}_+ \quad \forall k,\forall t \label{X}\\
&&P_k^{tdh},P_{k,\rm ch}^{tdh},P_{k,\rm dis}^{tdh},E_k^{tdh},P_{\rm lc}^{tdh} \in \mathbb{R}_+ \quad \forall k,\forall t, \forall d, \forall h \label{Y}
\end{eqnarray}

In the above formulation, the objective function \eqref{OBJ} minimizes the total budget ($TB$) of microgrid over the planning horizon. It consists of the investment cost ($C_{\rm inv}^t$) and operational cost ($C_{\rm opr}^t$). In \eqref{OBJ_1}, the investment costs, which are proportional to the nominal capacity and installation number of DERs, are calculated on the yearly basis. In \eqref{OBJ_2}, the operational costs include fixed operation and maintenance (O\&M) cost, fuel consumption cost, and load curtailment cost. The first term, O\&M cost, is estimated as the fixed yearly rate of DERs' installation. The rest terms are calculated on the daily basis and then scaled to yearly values. The fuel cost of DFGs is quantified using the linear cost factor that synthesized from the fuel consumption rate and efficiency~\cite{Moradi2015Operational} of each unit. The load curtailment cost is expressed as the product of load reduction and penalty factor. Overall, $TB$ is evaluated as the net present value, where $\tau(t)=(1+r)^{-t}$ denotes the present value factor.

Constraints \eqref{TB}-\eqref{QB2} are investment constraints. The annual investment cost is constrained by \eqref{TB}. Constraint \eqref{YearGap} imposes the logic relationship of DERs' installation over consecutive years. The number of installed DERs in each year is constrained by \eqref{XLimit}. Constraint \eqref{QB1} keeps the annual installation of DERs larger than the summation of peak load and long-term reserve. Constraint \eqref{QB2} regulates the minimum installation of DFGs, which serve as the stable power sources in a stand-alone system. Constraints \eqref{PB}-\eqref{Einitial} are operational constraints. Constraint \eqref{PB} ensures that the hourly power supply can fully cover the load demand in microgrid. Note that load curtailment ($P_{\rm lc}^{tdh}$) performs as a slack variable in \eqref{PB}, whose value can be naturally confined by imposing a large penalty factor ($q_{\rm v}$). In \eqref{Pres}, the generation of RESs is constrained by their available capacities, which are relevant to real-time resource conditions, e.g., wind speed and solar irradiation. The generic constraints of dispatchable generators are presented in \eqref{Pdg}-\eqref{Preserve}. The allowable generation range and hourly ramping capability of DFGs are constrained by \eqref{Pdg} and \eqref{Pramp}, respectively. Constraint \eqref{Preserve} ensures that the short-term reserve provided by DFGs can satisfy the operational requirement of microgrid~\cite{Zhao2013Short}. The operational states of ESS are constrained by \eqref{Pbes1}-\eqref{Einitial}~\cite{Hajipour2015Stochastic}. Since the operation of ESS is evaluated in independent days, their energy states in the first and last hours of each day should be kept the same. Finally, the decision variables in investment and operational problems are defined in \eqref{X} and \eqref{Y}, respectively.

As mentioned, MMEP problem is associated with a variety of random and non-random uncertainties. In \eqref{PB}, we assume that the load demand can be decomposed into demand capacity ($\psi^t$) and load variation factor ($p_{\rm L}^{tdh}$)~\cite{Rahmani2016Distributed}. The primary source of non-random uncertainties results from the forecast error of demand capacity growth, which can be described using the envelope bound model~\cite{Kazemi2015Risk}:
\begin{equation}
\label{EBM}
\Gamma(\alpha_{\rm L},\tilde{\psi})=\left\{\psi:|(\psi_t-\tilde{\psi}_t)/\tilde{\psi}_t| \leq \alpha_{\rm L} \ \mid \alpha_{\rm L}\geq 0, \forall t \right\}
\end{equation}
with $\tilde{\psi}_t$ being the forecasted value of demand capacity in each year. $\alpha_{\rm L}$ is the horizon of demand growth uncertainty, which defines the information gap between the forecasted demand levels and the actual values. Generally, for a fixed $\tilde{\psi}_t$, a greater $\alpha_{\rm L}$, i.e., a stronger risk-tolerance capability, requires a more capacitated and costly planning solution for the microgrid.

On the other hand, the load variation factor ($p_{\rm L}^{tdh}$) in each hour can be considered as the high-frequency component of load demand, which is classified as random uncertainties. Other types of random variables include the hourly power generation of RESs ($A_{\rm res}^{tdh}$), which refers to WT and PV in this paper. Detailed probabilistic models of these random factors can be found in~\cite{Billinton2008Effects,Ekren2008Size}.

Next, to better reflect the impact of random and non-random uncertainties, we extend our basic deterministic model to a chance constrained stochastic program.

\section{Formulation of Chance Constrained Information-Gap Decision Model}
\label{Sec:Formulation}
\subsection{Deterministic MMEP Model}
\label{Sec:DT}
\indent In the following, we first provide a matrix representation of the deterministic model assuming the nominal estimates of uncertain variables to avoid a cumbersome exposition.
\setlength{\arraycolsep}{0.0em}
\begin{eqnarray}
\min \ && \Lambda(x,y,\tilde{\psi})=c^{T}x+d^{T}y  \label{OBJ_dt}\\
s.t.\ && Ax\geq b  \label{Cons1_dt}\\
&&Hx-W\tilde{\psi}\geq h \label{Cons2_dt}\\
&&Sy-L\tilde{\psi}\geq 0 \label{Cons3_dt}\\
&&-Fx+Ty\geq 0  \label{Cons4_dt}\\
&&Gy\geq g  \label{Cons5_dt}\\
&&x\in \mathbb{Z}_+,y\in \mathbb{R}_+  \label{Cons6_dt}
\end{eqnarray}

\indent In \eqref{OBJ_dt}, the budget function of deterministic model is denoted by $\Lambda(x,y,\tilde{\psi})$, where $x$ and $y$ are the vectors of investment and operational variables, respectively. The variables $x$ are constrained by \eqref{Cons1_dt}, corresponding to constraints \eqref{TB}-\eqref{XLimit}. Constraints \eqref{Cons2_dt} and \eqref{Cons3_dt} correlate $x$ and $y$ with the forecasted value of load capacity growth ($\tilde{\psi}$), corresponding to constraints \eqref{QB1}-\eqref{QB2} and \eqref{PB} respectively. Constraint \eqref{Cons4_dt} builds the relationship between $x$ and $y$, corresponding to constraints \eqref{Pres}-\eqref{Ebes}. Rest constraints concerning $y$ are represented by \eqref{Cons5_dt}.

\begin{remark}
\label{rmk:0}
For an easy description, some equality constraints in \eqref{TB}-\eqref{Einitial} are represented as inequalities in \eqref{Cons1_dt}-\eqref{Cons5_dt}. Note that there is no loss of generality as one equality (e.g. $a=b$) can always be equivalently  replaced by a pair of inequalities (e.g. $a\leq b$ and $a\geq b$).
\end{remark}

\subsection{Information-Gap Decision Model for Microgrid Planning}
\label{Sec:IGD}
\indent In this paper, information gap decision theory (IGDT) is applied to address the non-random uncertainties in MMEP. The main idea of IGDT is to find a solution with maximized horizon of uncertainty as well as acceptable system performance~\cite{Kazemi2015Risk}. To hedge against the load uncertainty, the information-gap decision (IGD) model can be written as:
\setlength{\arraycolsep}{0.0em}
\begin{eqnarray}
\max\limits_{x,y} \ && \alpha_{\rm L}  \label{OBJ_igd}\\
s.t. \ && \bar{\Lambda} \leq (1+\sigma)\Lambda_0  \label{Cons1_igd}\\
&&Hx-W\psi\geq h   \label{Cons2_igd}\\
&&Sy-L\psi\geq 0 \label{Cons3_igd}\\
&&\rm Eqs. \ \eqref{Cons1_dt},\eqref{Cons4_dt}-\eqref{Cons6_dt}  \label{Cons4_igd}\\
&&\bar{\Lambda}=\max\limits_{\psi\in\Gamma(\alpha_{\rm L},\tilde{\psi})}\left\{\Lambda(x,y,\psi)\right\} \label{Cons5_igd}
\end{eqnarray}

\indent {{Note that the IGD model is essentially a bi-level optimization problem.}} The upper level in \eqref{OBJ_igd} aims to maximize the horizon of load uncertainty meanwhile satisfying constraints \eqref{Cons1_igd}-\eqref{Cons4_igd}. Constraint \eqref{Cons1_igd} is named as the 'budget level limit', which keeps the total budget ($\bar{\Lambda}$) lower than a specified level $(1+\sigma)\Lambda_0$. In this expression, $\Lambda_0$ represents the risk-neutral budget level, while $\sigma$ is the deviation factor that specifies the acceptable degree of budget excess. Constraints \eqref{Cons2_igd} and \eqref{Cons3_igd} associate the decision variables with demand growth uncertainty ($\psi$). In \eqref{Cons5_igd}, the lower level determines the robust performance of MMEP, i.e., the required budget level, to hedge against any possible variation of $\psi$ in the set of $\Gamma(\alpha_{\rm L},\tilde{\psi})$.

\indent {{To address the computational challenge from bi-level structure, we convert it into a tractable single-level form by fixing the demand growth at its maximum value, i.e., $\psi=(1+\alpha_{\rm L})\tilde{\psi}$. This conversion is typically valid according to its deterministic counterpart in Section II, where the largest budget requirement always occurs when demand growth achieves its upper bound.}} Hence, the bi-level structure can be simplified by replacing \eqref{Cons1_igd}-\eqref{Cons3_igd} and \eqref{Cons5_igd} with following \eqref{Cons6_igd}-\eqref{Cons8_igd}.
\setlength{\arraycolsep}{0.0em}
\begin{eqnarray}
&&c^{T}x+d^{T}y \leq (1+\sigma)\Lambda_0  \label{Cons6_igd}\\
&&Hx-(1+\alpha_{\rm L})W\tilde{\psi}\geq h  \label{Cons7_igd}\\
&&Sy-(1+\alpha_{\rm L})L\tilde{\psi}\geq 0 \label{Cons8_igd}
\end{eqnarray}


{{As implied by \eqref{OBJ_igd}, \eqref{Cons4_igd}, \eqref{Cons6_igd}-\eqref{Cons8_igd}, the solution of IGD model could be evidently influenced by the value of deviation factor. There is a possibility of over-conservatism due to the inappropriate selection of model parameter. In addition, the IGDT-based formulation becomes less attractive when more than one types of uncertainties are involved, which may require for a multi-objective framework~\cite{Dehghan2014Multi}. To overcome these shortcomings, we develop a scenario-based chance constrained extension of IGD model in the next section.}}

\subsection{Chance Constrained Information-Gap Decision Model}
\label{Sec:CCIGD}
We mention that IGD model does not consider the random factors, which, however, is very critical in the operational stage. So, to further consider the random uncertainties ($\xi$) in MMEP, we replace the coefficient matrices in \eqref{Cons4_dt} and \eqref{Cons8_igd} by a set of stochastic matrices, corresponding to a set of stochastic scenarios. Definitely, the budget level limit \eqref{Cons1_igd} and power supply constraints \eqref{Cons8_igd} can be imposed for every realization of $\xi$. However, this requirement could be too restrictive, and practically they are allowed to be violated under extreme situations.  Hence, we impose a chance constraint to ensure that \eqref{Cons1_igd} and \eqref{Cons8_igd} are satisfied with a predefined probability. By doing so, the essential part of chance constrained IGD model can be written as:
\setlength{\arraycolsep}{0.0em}
\begin{eqnarray}
&&\mathbb{P}\left\{
\begin{array}{c}
c^{T}x+d^{T}y(\xi)\leq (1+\sigma)\Lambda_0 \\
Sy(\xi)-(1+\alpha_{\rm L})L(\xi)\tilde{\psi}\geq 0
\end{array}
\right\} \geq 1-\varepsilon \label{RC-a}\\
&&-F(\xi)x+Ty(\xi)\geq 0 \label{RC-c} \\
&&Gy(\xi)\geq g \label{RC-d}
\end{eqnarray}

Note that the chance constrained IGD model takes a step forward by incorporating random $\xi$ in the operational stage,
whose joint probability space is defined in set $\Omega\subseteq\mathbb{R}$. Constraint \eqref{RC-a} requires that the budget level limit and power supply constraints should be satisfied with a probability no less than $1-\varepsilon$, where $\varepsilon$ denotes the risk tolerance level. In \eqref{RC-a} and \eqref{RC-c}, random matrices $L(\xi)$ and $F(\xi)$ are associated with the distributions of load variation and RES generation, respectively.

To solve this type of problem, one popular approach is to convert the chance constrained model into a deterministic equivalence if the $1-\varepsilon$ quantile of the underlying distribution can be analytically represented. However, this method is not applicable for the joint probabilistic constraints, which demand for multi-variate integration. Indeed, our chance constraint \eqref{RC-a} combines a set of constraints and several random variables, which  makes it infeasible to derive a closed-form expression of the $1-\varepsilon$ quantile. Thus, we adopt a sampling strategy to represent the randomness by  a finite set of discrete scenarios $\Omega=\{\xi_1,\xi_2,...,\xi_N\}$, and  formulate  a novel bilinear model of chance constrained program~\cite{Zeng2014Chance}. The bilinear formulation of chance constrained IGD model can be written as:
\setlength{\arraycolsep}{0.0em}
\begin{eqnarray}
\max \ && \alpha_{\rm L} \label{OBJ_ccigd}\\
s.t. \ && \rm Eqs. \eqref{Cons1_dt},\eqref{Cons7_igd} \label{Cons0_ccigd}\\
&&[c^{T}x+d^{T}y_n-(1+\sigma)\Lambda_0](1-z_n)\leq 0 \quad \forall n \label{Cons1_ccigd}\\
&&[Sy_n-(1+\alpha_{\rm L})L_n\tilde{\psi}](1-z_n)\geq 0 \quad \forall n \label{Cons2_ccigd}\\
&&-F_n x+Ty_n\geq 0 \quad \forall n  \label{Cons3_ccigd}\\
&&Gy_n\geq g \quad \forall n \label{Cons4_ccigd}\\
&&\sum_{n=1}^N \pi_n z_n\leq\varepsilon \label{Cons5_ccigd}\\
&&\alpha_{\rm L}\geq 0,x\in\mathbb{Z}_+ \label{Cons6_ccigd}\\
&&y_n\in\mathbb{R}_+,z_n\in\{0,1\},n=1,...,N \label{Cons7_ccigd}
\end{eqnarray}

In the aforementioned model, $\alpha_{\rm L}$ and $x$ are first-stage variables, while $y_n$ denotes second-stage operations (i.e., recourse) variables. The budget level limit is decomposed into a set of constraints. All the budget level constraints share the same first-stage cost $c^{T}x$, but with separated recourse cost $d^{T}y_n$ defined in each scenario. The first-stage problem is only related to non-random uncertainties, which is subject to constraints \eqref{Cons0_ccigd}. The second-stage recourse problem is associated with both random and non-random uncertainties, which is subject to constraints \eqref{Cons1_ccigd}-\eqref{Cons5_ccigd}.

Chance constraint \eqref{RC-a} is replaced by the bilinear structure in \eqref{Cons1_ccigd} and \eqref{Cons2_ccigd}. The binary variable $z_n$ is employed to indicate the full requirements of scenario $n$ is imposed or not. When $z_n=1$, \eqref{Cons1_ccigd} and \eqref{Cons2_ccigd} can be ignored, i.e., those constraints are deactivated, and only \eqref{Cons3_ccigd} and \eqref{Cons4_ccigd} are necessary. Otherwise, the complete constraints of scenario $n$, i.e., \eqref{Cons1_ccigd}-\eqref{Cons4_ccigd}, must be satisfied. The total probability of scenarios that are partially satisfied is restricted by \eqref{Cons5_ccigd}, where $\pi_n$ denotes the realization probability of scenario $n$.


\begin{remark}
\label{rmk:1}
Note that if $z_n=1$, which means that constraints in \eqref{Cons1_ccigd} and \eqref{Cons2_ccigd} will be ignored. For the remaining constraints associated with scenario $n$, i.e., \eqref{Cons3_ccigd} and \eqref{Cons4_ccigd}, we mention that it is always feasible to define a feasible $y_n$ for any given $x$. The reason behind is that \eqref{Cons3_ccigd} and \eqref{Cons4_ccigd}, which represent \eqref{Pres}-\eqref{Einitial} in the basic deterministic model, define the feasible set for every single operation variable that is never empty. Thus, we can easily generate $y_n$ satisfying those constraints for any fixed $x$. Moreover, given that the objective function in \eqref{OBJ_ccigd} does not depend on $y_n$, we can conclude that \eqref{Cons3_ccigd} and \eqref{Cons4_ccigd} alone do not impact the choice of the first stage decision variables ($\alpha_{\rm L},x$). Note that such observation allows us to simply focus on the impact of \eqref{Cons1_ccigd} and \eqref{Cons2_ccigd} by making decisions on the selection of scenario $n$.
\end{remark}

Overall, the CC-IGD based microgrid planning model in \eqref{OBJ_ccigd}-\eqref{Cons7_ccigd} gives a full consideration to the random and non-random uncertainties defined in Section II. On one hand, the multi-period investment decisions are made to maximize the risk-tolerance towards long-term uncertainty of load growth. On the other hand, the hourly operational decisions are made under each scenario to hedge against the random uncertainties of RES generation and load variation. Next, we present an effective computational method to solve this sophisticated model.

\section{Solution Approach}
\indent The mixed-integer nonlinear formulation presented in \eqref{OBJ_ccigd}-\eqref{Cons7_ccigd} for MMEP could be very difficult to compute, even for the case with a small number of scenarios to capture randomness. To address such computational challenge, we follow the idea presented in~\cite{Zeng2014Chance} to develop a strengthened bilinear Benders decomposition algorithm to solve the proposed planning problem. As a decomposition method, our algorithm involves a master problem and two subproblems. Before proceed to concrete algorithm steps, we describe subproblems in the follows.

Our first subproblem, which is referred to the economic dispatch subproblem (EDS) for a single scenario, is to minimize the recourse cost $d^{T}y_n$ subject to other constraints in the second stage, i.e., \eqref{Cons2_ccigd}-\eqref{Cons4_ccigd}, for a given first stage solution, denoted by $(\hat{\alpha}_{\rm L},\hat{x})$. Mathematically, it is
\setlength{\arraycolsep}{0.0em}
\begin{eqnarray}
\mathbf{SP_1}:    \  \hat{\Phi}_n=&& \min \ d^{T}y_n \label{OBJ_osp}\\
s.t. \ &&  Sy_n \geq(1+\hat{\alpha}_{\rm L})L_n\tilde{\psi}:\mu_n  \label{Cons1_osp}\\
&&Ty_n\geq F_n \hat{x}:\lambda_n   \label{Cons2_osp}\\
&&Gy_n\geq g:\gamma_n              \label{Cons3_osp}\\
&&y_n,\mu_n,\lambda_n,\gamma_n\in\mathbb{R}_+  \label{Cons4_osp}
\end{eqnarray}
where $\mu_n,\lambda_n,\gamma_n$ are dual variables of constraints \eqref{Cons1_osp}-\eqref{Cons3_osp}, respectively. As mentioned, the aforementioned
subproblem is defined for every single scenario.

Note that EDS is always feasible due to load reduction variables $(P_{\rm lc}^{tdh})$ . According to the strategy in~\cite{Zeng2014Chance},  we just need to supply a Benders optimality cut ($OC$) in the following bilinear form to the master problem.
\begin{equation}
\label{eq_OC}
[\hat{\Phi}_n+\hat{\mu}_n^{T}L_n(\alpha_{\rm L}-\hat{\alpha}_{\rm L})+\hat{\lambda}_n^{T}F_n(x-\hat{x})](1-z_n)\leq\Phi_n
\end{equation}
where variable $\Phi_n$ is to represent EDS cost in the second stage and $z_n$ is the binary indicator variable defined in \eqref{OBJ_ccigd}-\eqref{Cons7_ccigd}.

\begin{remark}
\label{rmk:2}
$(i)$ Unlike the traditional Bender decomposition method, $OC$ is modulated by the binary indicator variable $z_n$. When $z_n=1$, \eqref{eq_OC} reduces to a trivial constraint $\Phi_n\geq 0$, which means $OC$ is deactivated. Otherwise, it is enforced in the master problem. \\
$(ii)$ As mentioned earlier, when $z_n=1$, only \eqref{Cons3_ccigd} and \eqref{Cons4_ccigd}  are presented for scenario $n$. As they do not affect either the feasibility or the optimality of $(\hat{\alpha}_{\rm L},\hat{x})$, we are not concerned by those constraints, except for their impact in $\mathbf{SP_1}$.
\end{remark}

Another subproblem is designed to verify the feasibility of the first stage solution $(\hat{\alpha}_{\rm L},\hat{x})$ to the original formulation in \eqref{OBJ_ccigd}-\eqref{Cons7_ccigd}. Let $t_n$ be the $n$th auxiliary variable. We have
\setlength{\arraycolsep}{0.0em}
\begin{eqnarray}
\mathbf{SP_2}: \ \Delta = && \min \ \sum_{n=1}^N t_n  \label{OBJ_ogp}\\
s.t. \ && t_n\geq [c^{T}\hat{x}+\hat{\Phi}_n-(1+\sigma)\Lambda_0](1-z_n) \ \forall n   \label{Cons2_ogp}\\
&&\sum_{n=1}^N \pi_n z_n \leq \varepsilon \label{Cons3_ogp}\\
&&t_n\geq 0, z_n\in\{0,1\}, \ n=1,2,...,N. \label{Cons4_ogp}
\end{eqnarray}

\begin{remark}
\label{rmk:3}
We can easily see that once $\Delta=0$, $(\hat{\alpha}_{\rm L},\hat{x})$ is a feasible solution of chance constrained IGD model. This property is actually critical to our Benders decomposition algorithm development. As we will show next, if an optimal solution of the master problem leads to zero optimal value in $\mathbf{SP_2}$, it is also optimal to the chance constrained IGD model.
\end{remark}

In the following, we define the master problem ($\mathbf{IGDMP}$) based on the strategy provided in~\cite{Zeng2014Chance}.  Let $U_{n,i}(\alpha_{\rm L},x)=\hat{\Phi}_{n,i}+\hat{\mu}_{n,i}^{T}L_n(\alpha_{\rm L}-\hat{\alpha}_{{\rm L},i})+\hat{\lambda}_{n,i}^{T}F_n(x-\hat{x_i})$, where $i$ is the counter of iteration.
\setlength{\arraycolsep}{0.0em}
\begin{eqnarray}
\max \  && \alpha_{\rm L}  \label{OBJ_igdmp}\\
s.t. \  && Ax\geq b  \label{Cons1_igdmp}\\
&&Hx-(1+\alpha_{\rm L})W\tilde{\psi}\geq h  \label{Cons2_igdmp}\\
&&c^{T}x+\Phi_n-(1+\sigma)\Lambda_0\leq 0 \quad \forall n \label{Cons3_igdmp}\\
&&(1-z_n)U_{n,i}(\alpha_{\rm L},x)\leq\Phi_n \ i=1,...,j-1,\forall n \label{Cons4_igdmp}\\
&&\sum_{n=1}^N \pi_n z_n \leq \varepsilon \label{Cons6_igdmp}\\
&&\alpha_{\rm L}\geq 0,x\in \mathbb{Z}_+ \label{Cons7_igdmp}\\
&&z_n\in\{0,1\},\Phi_n\in\mathbb{R}_+,n=1,...,N \label{Cons8_igdmp}
\end{eqnarray}

Similar to the conventional Benders decomposition method, $OC$ in \eqref{Cons4_igdmp}, in a bilinear form, will be generated from $\mathbf{SP}_1$ for every scenario and  then  supplied to this master problem in every iteration.

\begin{remark}
\noindent $(i)$ The bilinear structure in \eqref{Cons4_igdmp} can be easily linearized by using McCormick linearization method~\cite{Mccormick1976Computability}. For example, consider the bilinear term $x_k z_n(=\tilde{x}_k)$, where $x_k$ is the $k$-th component variable of $x$. It can be replaced by the following linear constraints and variables, given that $\bar{x}_k$ is an upper bound of $x_k$. As a result, the master problem is converted into an MIP problem that can readily be computed by any professional MIP solver.
\begin{equation}
\label{mcl}
\tilde{x}_k\leq \bar{x}_k z_n,\ \tilde{x}_k\leq x_k,\ \tilde{x}_k\geq x_k-\bar{x}_k(1-z_n),\ \tilde{x}_k\geq 0
\end{equation}
\noindent $(ii)$ Nevertheless, the linearization in \eqref{mcl} will introduce a large amount of new variables and constraints, which makes it computational cumbersome to solve the $\mathbf{IGDMP}$. To address this shortcoming, we adopt an equivalent replacement for the bilinear formulation in \eqref{Cons3_igdmp} and \eqref{Cons4_igdmp}:
\setlength{\arraycolsep}{0.0em}
\begin{eqnarray}
&&c^{T}x+\Phi_n(1-z_n)-(1+\sigma)\Lambda_0\leq 0 \quad \forall n \label{Cons1_nv}\\
&&U_{n,i}(\alpha_{\rm L},x)\leq\Phi_n \quad i=1,...,j-1,\forall n \label{Cons2_nv}
\end{eqnarray}
The equivalence can be argued by considering \eqref{Cons3_igdmp}-\eqref{Cons4_igdmp} and \eqref{Cons1_nv}-\eqref{Cons2_nv} in both $z_n=0$ and $z_n=1$ situations. When $z_n=0$, \eqref{Cons3_igdmp} and \eqref{Cons4_igdmp} reduce to \eqref{Cons1_nv} and \eqref{Cons2_nv}. And when $z_n=1$, $\Phi_n\geq 0$ and hence $\Phi_n$ can be simply set to 0 in \eqref{Cons3_igdmp},  which indicates that an optimal $(\alpha_L,x)$ to $\mathbf{IGDMP}$ is not affected by $U_{n,i}(\alpha_{\rm L},x)$ or $\Phi_n$. Note that the same effect is also achieved in \eqref{Cons1_nv} due to $\Phi_n(1-z_n)=0$.

We mention that linearization of \eqref{Cons1_nv} requires a much less number of newly-added variables and constraints than that from linearization of \eqref{Cons4_igdmp}. Hence, the modified $\mathbf{IGDMP}$ with \eqref{Cons1_nv}-\eqref{Cons2_nv} could be solved in a more efficient way. In the rest of this paper, we denote Benders decomposition with this improvement as \textbf{strengthened bilinear Benders Decomposition} while the one with \eqref{Cons3_igdmp} and \eqref{Cons4_igdmp} as the original one.

\noindent $(iii)$ Note that the master problem is a relaxation to the chance constrained IGD model in \eqref{OBJ_ccigd}-\eqref{Cons7_ccigd}, which will be strengthened by iteratively adding  \eqref{Cons4_igdmp}.  Consider its optimal solution $(\hat{\alpha}_{\rm L},\hat{x})$. If the corresponding $\mathbf{SP}_2$ returns zero, which verifies that it is feasible to \eqref{OBJ_ccigd}-\eqref{Cons7_ccigd}, we can conclude that it is also optimal to \eqref{OBJ_ccigd}-\eqref{Cons7_ccigd}.
\end{remark}

With the well-defined subproblems and the master problem, the complete steps of our method is outlined in \textbf{Algorithm} \ref{alg:mbbd}.

\begin{algorithm}[!t]
\caption{Strengthened Bilinear Benders Decomposition}
\begin{algorithmic}[1]
\REQUIRE ~~\\
Set $j=0$ and $\Delta_j=+\infty$;$OC_{n,j}=\emptyset$;
\ENSURE ~~\\
\WHILE    {$\Delta_j>0$}
  \STATE $\mathbf{IGDMP}_j\leftarrow\mathbf{IGDMP}_{j-1} \bigcup\limits_{n=1}^N \{OC_{n,j}\}$;
  \STATE Compute master problem $\mathbf{IGDMP}_j$;
     \IF {$\mathbf{IGDMP}_j$ is infeasible}
     \STATE terminate and \emph{report the infeasibility};
     \ELSE
     \STATE Get solution $(\hat{\alpha}_{\rm L}^j,\hat{x}^j,\hat{z}^j)$;
     \ENDIF
  \STATE $j\leftarrow j+1$;
     \FOR {$n=1$ to $N$}
     \STATE  Compute $\mathbf{SP}_1$ to get $(\hat{\Phi}_{n,j},\hat{\mu}_{n,j},\hat{\lambda}_{n,j})$;
     \STATE  Generate optimality cut $OC_{n,j}$;
     \ENDFOR
   \STATE Compute $\mathbf{SP}_2$ to get $\Delta_j$;
\ENDWHILE
\STATE {\bf output:} Report $(\hat{\alpha}_{\rm L}^j,\hat{x}^j)$.
\end{algorithmic}
\label{alg:mbbd}%
\end{algorithm}

\begin{remark}
\noindent
The algorithm of Strengthened Bilinear Benders Decomposition terminates with an optimal solution to \eqref{OBJ_ccigd}-\eqref{Cons7_ccigd} or report infeasibility in a finite number of iterations.
\end{remark}

Note that whenever $\mathbf{SP}_2$ has a strictly positive optimal value, i.e., the latest optimal solution to the master problem is not feasible to the original CC-IGD model, a set of new extreme points (i.e., their corresponding OC cuts) will be included into the master problem $\mathbf{IGDMP}$. It must differ from those generated in previous iterations. Otherwise, the new master problem is same as that in the previous iteration. However, given the connection between \eqref{Cons3_igdmp}-\eqref{Cons4_igdmp} and $\mathbf{SP}_2$, we can conclude that both master problems (because they are identical) are infeasible, which, according to \textbf{Algorithm} \ref{alg:mbbd}'s description, terminates the whole algorithm.

Note also that the feasible set of dual problem to $\mathbf{SP}_1$ has a finite number of extreme points, which is the basis of Benders $OC$ cuts, for every scenario. Hence, we conclude that the algorithm will either produce an optimal solution or  terminate with infeasibility  in a finite number of iterations.

%
%
\section{Numerical Results}
In this section, the proposed planning model and computation approach are firstly demonstrated in a simple microgrid through a series of comparative numerical tests. Then, we further test our method on a practical microgrid with complex structure to prove its scalability. All the algorithm development and computations, including our bilinear Benders decomposition, are made by CPLEX in MATLAB environment on a desktop computer with Intel Core i5-3470 3.20 GHZ processer.

\subsection{Simple Microgrid}
\label{Sec:Case1}
In this case, the proposed method is tested on a stand-alone microgrid with a simple structure. The candidate DERs include wind turbines (WT), photovoltaic panels (PV), diesel engine (DE), micro turbine (MT), and battery packs (ES), whose major parameters are listed in Table \ref{tab:der1}. The planning horizon is 10 years with 5 separated periods. The discount rate is set to 4\%. The initial demand capacity is 450kW and its growth rate is forecasted to be 8\% per year. The penalty cost for load curtailment is specified as 5 \$/kWh. Four typical days are selected from each year to represent the seasonal and diurnal patterns of RES generation and load variation. According to the historical data in~\cite{Cao2016Multi}, 2000 scenarios are produced using Monte Carlo simulation. Overall, each scenario contains three time-series for wind, solar, and load respectively, each of which has totally 960 sampling points.

\renewcommand\arraystretch{0.7}
\begin{table}[!t]
  \centering
  \caption{Simple Microgrid: Parameters of Candidate DERs}
  \setlength{\tabcolsep}{0.30mm}
    \begin{tabular}{ccccccccc}
    \toprule
    \multicolumn{9}{c}{Distributed Generators} \\
    \midrule
    \multirow{3}[2]{*}{Label} & \multirow{3}[2]{*}{Type} & Rated & Min & Capital & O\&M  & Fuel  & Num &  \\
          &       & Power & Output & Cost  & Cost  & Cost  & Limit &  \\
          &       & (kW)  & (kW)  & (\$/kW) & (\$/kW/yr) & (\$/kWh) &       &  \\
    \midrule
    WT    & Wind Turbine & 120   & 0     & 900   & 0     & /     & 3     &  \\
    \midrule
    PV    & Photovolatic & 80    & 0     & 1200  & 0     & /     & 3     &  \\
    \midrule
    DE    & Diesel Engine & 60    & 10    & 300   & 175.2 & 0.105 & 2     &  \\
    \midrule
    MT    & Micro Turbine & 80    & 10    & 450   & 262.8 & 0.089 & 2     &  \\
    \midrule
    \multicolumn{9}{c}{Energy Storage Devices} \\
    \midrule
    \multirow{3}[2]{*}{Label} & \multirow{3}[2]{*}{Type} & Rated & Rated & Power & Energy  & O\&M  & Effi- & Num \\
          &       &  Power & Energy & Cost  & Cost  & Cost  & \multirow{2}[1]{*}{ciency} & \multirow{2}[1]{*}{Limit} \\
          &       & (kW)  & (kWh) & (\$/kW) & (\$/kWh) & (\$/kW/yr) &       &  \\
    \midrule
    ES1   & Battery Pack 1 & 90    & 150   & 270   & 150   & 23.4  & 0.95  & 4 \\
    \midrule
    ES2   & Battery Pack 2 & 100   & 200   & 200   & 180   & 35    & 0.85  & 4 \\
    \bottomrule
    \end{tabular}%
  \label{tab:der1}%
\end{table}%

First, we will exhibit the basic planning results. To make a compromise between the computational accuracy and efficiency, scenario reduction is performed to extract 80 representative scenarios for problem solving. Next, the risk-averse capability of proposed CC-IGD model will be compared with deterministic and IGD-based MMEP models through performance evaluation, where the entire scenario set is included. Moreover, sensitivity analysis will be performed to investigate the impact of key control parameters, i.e., deviation factor $\sigma$ and risk tolerance level $\varepsilon$, on planning results. Last but not the least, the computational feasibility of strengthened bilinear Benders decomposition method will be tested.

\subsubsection{Basic Planning Results}
To better illustrate the planning results of CC-IGD, we choose deterministic MMEP model (DT) in Section \ref{Sec:DT} and IGD model (IGD) in Section \ref{Sec:IGD} as benchmarks. The parameters are set as deviation factor $\sigma=0.4$, risk tolerance level $\varepsilon=0.1$. Under given conditions, the objective value of DT is minimized as $\Lambda_0=\$5,425,087$, which provides the risk-neutral budget for IGD and CC-IGD computation. Set the budget limit as $1.4\times \$5,425,087$, the robustness level ($\alpha_{\rm L}$) maximized by IGD and CC-IGD are 0.1742 and 0.0857, respectively. CC-IGD derives lower $\alpha_{\rm L}$ due to its extra consideration for random factors. Table \ref{tab:sd} presents the solutions of different MMEP models. We regard DT solution as the basic scheme, while the solutions of IGD and CC-IGD as the reinforcement. To prepare for the unexpected demand growth, IGD makes comprehensive reinforcement by installing 200kW more RES, 100kW more DFGs, and 200kWh more ESS. Different from the IGD solution, the reinforcement made by CC-IGD seems to be more selective. First, CC-IGD brings forward the installation of 60kW DE from Period V to Period II. Then, the WT installation is reduced from 600kW to 240kW at Period II \& III. Instead, CC-IGD adds the capacity of MT from 80kW to 160kW at Period III. All these changes could reduce the random risks brought by RES output and load variation while retaining system's robustness against non-random uncertainties.
\renewcommand\arraystretch{0.7}
\begin{table*}[!t]
  \centering
  \caption{Simple Microgrid: Planning Results}
  \setlength{\tabcolsep}{3.0mm}
    \begin{tabular}{ccccc}
    \toprule
    Period & DER Type & \#CC-IGD\# & DT    & IGD \\
    \midrule
    \multirow{3}[2]{*}{I
Year 1-2} & RES   & WT=360kW PV=240kW & WT=360kW PV=240kW & WT=360kW PV=240kW \\
          & DFG   & DE=120kW MT=160kW & DE=120kW MT=160kW & DE=120kW MT=160kW \\
          & BES   & ES1=600kWh ES2=200kWh & ES1=600kWh ES2=400kWh & ES1=600kWh ES2=200kWh \\
    \midrule
    \multirow{3}[2]{*}{II
Year 3-4} & RES   & WT=240kW PV=240kW & WT=360kW PV=240kW & WT=360kW PV=240kW \\
          & DFG   & DE=60kW MT=160kW & MT=160kW & MT=160kW \\
          & BES   & ES1=450kWh ES2=600kWh & ES1=600kWh ES2=800kWh & ES1=600kWh ES2=800kWh \\
    \midrule
    \multirow{3}[2]{*}{III
Year 5-6} & RES   & PV=240kW & WT=120kW PV=160kW & WT=240kW PV=240kW \\
          & DFG   & MT=160kW & 0     & MT=80kW \\
          & BES   & /     & ES2=400kWh & ES1=150kWh ES2=800kWh \\
    \midrule
    \multirow{3}[2]{*}{IV
Year 7-8} & RES   & WT=120kW & WT=120kW & WT=120kW \\
          & DFG   & DE=80kW & DE=120kW & MT=80kW \\
          & BES   & ES1=150kWh & ES1=150kWh & / \\
    \midrule
    V
Year 9-10 & DFG   & /     & /     & DE=60kW \\
    \bottomrule
    \end{tabular}%
  \label{tab:sd}%
\end{table*}%

\subsubsection{Performance Evaluation of Chance Constrained IGD Model}
The risk-averse capability of CC-IGD solution is further examined through post performance evaluation. The performance metric for planning schemes is defined as the Expected Project Budget (EPB), which can be calculated as:
\begin{equation}
EPB=c^{T}x+\frac{1}{NS}\sum_{n=1}^{NS} d^{T}y_n
\end{equation}
where variables $x$ and $y_n$ are defined as in Section III. $NS$ denotes the number of scenarios, which is set to 2000 in this case.

The performance metrics of CC-IGD solution are computed and compared with those of DT and IGD solutions under different robustness levels ($\alpha_{\rm L}$), as shown in Table \ref{tab:rbt}. The maximum load capacity during the planning horizon is intuitively presented in Column ``Peak''. When $\alpha_{\rm L}=0$ (Peak=900kW), the evaluation results of different models are close to each other. When $\alpha_{\rm L}=0.1$ (Peak=990kW), the $EPB$ of CC-IGD solution is \$7,065,800, which is below the budget level limit, i.e., \$1.4$\times$ 5,425,087=7,595,122. In contrast, the $EPB$ of DT solution exceeds the budget limit by \$322,878 due to the costly penalty for load curtailment. When $\alpha_{\rm L}=0.2$ (Peak=1,080kW), both DT and CC-IGD solutions violate the budget level limit. That is because $\alpha_{\rm L}=0.2$ is far beyond the maximum robustness level that can be tolerated by CC-IGD solution, i.e., 0.0857. However, the CC-IGD solution costs 20.70\% less than the DT. These observations prove the validity of the proposed CC-IGD model. Moreover, CC-IGD shows better performance than the IGD model in every instance, which reduces the budget expectation by more than \$195,800. Hence, CC-IGD shows a superiority over DT and IGD under the co-existence of random and non-random uncertainties.

\renewcommand\arraystretch{0.7}
\begin{table}[!t]
  \centering
  \caption{Simple Microgrid: Results of Performance Evaluation }
  \setlength{\tabcolsep}{3.0mm}
    \begin{tabular}{ccccc}
    \toprule
    \multirow{2}[4]{*}{$\alpha_{\rm L}$} & \multirow{2}[4]{*}{Peak} & \multicolumn{3}{c}{Expected Project Budget} \\
\cmidrule{3-5}          &       & DT    & IGD   & CC-IGD \\
    \midrule
    0.00  & 900kW & \$5,928,600 & \$6,043,800 & \$5,780,000 \\
    \midrule
    0.10  & 990kW & \$7,918,000 & \$7,261,600 & \$7,065,800 \\
    \midrule
    0.20  & 1080kW & \$10,904,000 & \$9,018,400 & \$8,647,100 \\
    \bottomrule
    \end{tabular}%
  \label{tab:rbt}%
\end{table}%

\subsubsection{Sensitivity Analysis of Key Control Parameters}
To analyze the impact of key control parameters, we solve the proposed planning model under different values of deviation factor ($\sigma$) and risk tolerance level ($\varepsilon$). Fig. 1(a) exhibits the varying trend of robustness level ($\alpha_{\rm L}$) of CC-IGD solution with respect to $\sigma$, which is the acceptable budget excess with respect to the risk-neutral budget from deterministic model. Fixing $\varepsilon=0.10$ while changing $\sigma$ from 0.20 to 0.50 with a step of 0.05, the robustness level improves from 0.0171 to 0.1314, which represents an increased immunity against uncertainties. This observation confirms an intuition that a larger $\sigma$ leads to a more capable planning solution to handle the long-term demand growth uncertainty.

Similarly, Fig. 1(b) presents the trend of robustness level via changing $\varepsilon$. As implied by CCP, a larger $\varepsilon$ leads a softer budget limit, i.e., more violations of budget limit are allowed over stochastic scenarios. Equivalently, more operational scenarios can be ignored in our MMEP planning model, which allows us to allocate more resources to maximize our capacity on handling the long-term demand growth uncertainty. Such insight is verified in our numerical study.  Fixing $\sigma=0.4$ while changing $\varepsilon$ from 0.00 to 0.30 with a step of 0.05, the robustness level grows from 0.0483 to 0.1462. Clearly, we can adjust $\varepsilon$ in CC-IGD model to achieve a balanced microgrid investment plan under both the long-term and operational uncertainties.
%
\begin{figure}[!t]
\centering
\subfloat[Deviation Factor]{\includegraphics[width=2.8in]{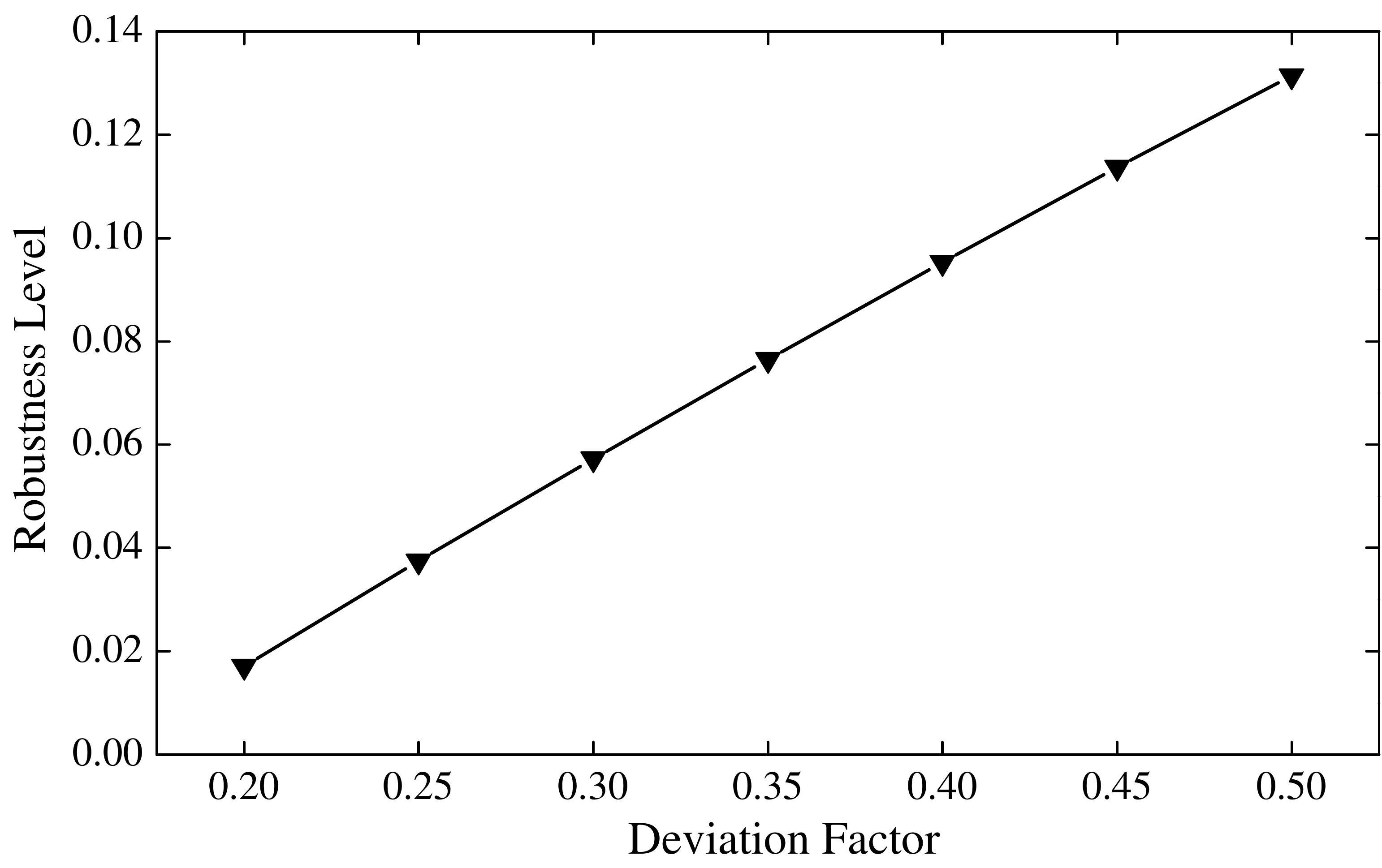}%
\label{fig:SDF}}
\hfil
\subfloat[Risk Tolerance Level]{\includegraphics[width=2.8in]{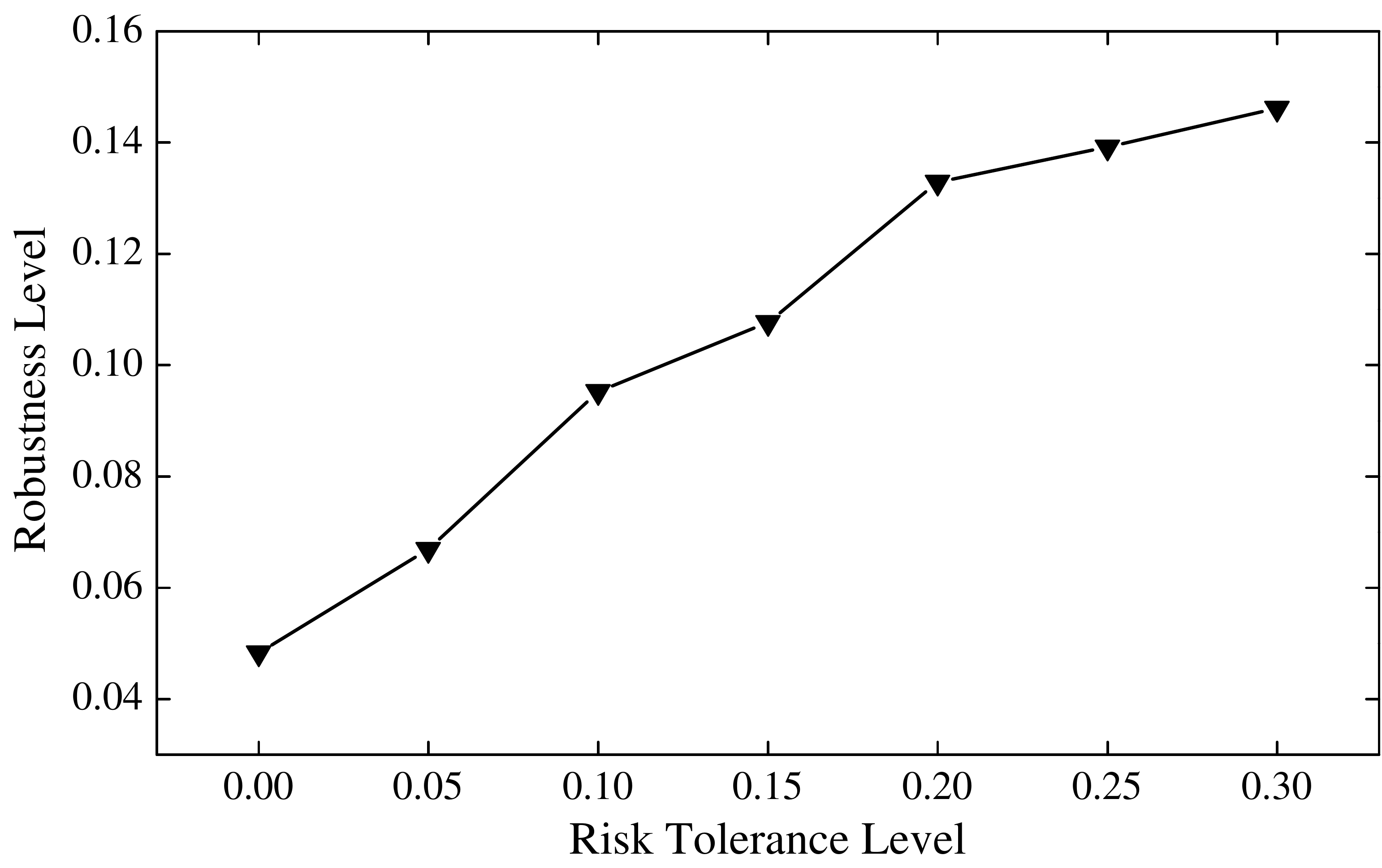}%
\label{fig:SRL}}
\caption{Sensitivity Analysis of Model Parameters}
\label{fig:SA}
\end{figure}

\subsubsection{Computational Test of Strengthened Bilinear Benders Decomposition Method}
The computational feasibility of SBD is verified by comparing to the original Benders decomposition (OBD) form, which adopts \eqref{Cons3_igdmp} and \eqref{Cons4_igdmp} in its master problem,  and the direct use of CPLEX (CPX). Again, McCormick linearization method is applied to convert \eqref{OBJ_ccigd}-\eqref{Cons7_ccigd} into mixed-integer linear formulation so that it can be directly solved by CPLEX. The solution time is counted by minutes, which has a limit up to 600 mins. If any problem is terminated due to the time limit, its solution time will be labeled as ``T''. The remaining gap or ``N/A'' (in case the gap cannot be reported) will be marked beside it.

Table \ref{tab:cp} presents the test results under different numbers of scenario ($NS$). The objective level, iteration number, solution time, and remaining gap are recorded in columns ``$\alpha_{\rm L}$'', ``itr'', ``min'', and ``gap'', respectively. The result indicates that bilinear BD method, even in its original form, performs better than the direct use of CPLEX. Only when $NS=20$ can all the algorithms converge to an optimal solution. SBD derives the same result as CPX but performs an order of magnitude faster, saving the solution time by 95.62\%. For the rest instances when scenario size gets larger, CPX either fails to derive any feasible solutions or cannot solve the linear programming relaxation within the time limit, which leads to no gap information available.

Also, the strengthened bilinear BD significantly outperforms its original counterpart (OBD). When $NS=20,40$, both methods succeed in solving the problem, but SBD is observed to be much faster. When $NS=60$ and $80$, only SBD can derive the optimal solution within the time limit. Compared to OBD, SBD might need more iterations to converge but achieves a tremendous time-saving in each iteration. Therefore, our modification in the bilinear Benders decomposition method has clearly improved its computational efficiency.

\renewcommand\arraystretch{0.7}
\begin{table}[!t]
  \centering
  \setlength{\tabcolsep}{1.4mm}
  \caption{Simple Microgrid: Solution Performance under Different Scenario Numbers}
    \begin{tabular}{cccccccccccc}
    \toprule
    \multirow{2}[4]{*}{\textit{NS}} & \multirow{2}[4]{*}{OBJ} & \multicolumn{3}{c}{\#SBD\#} & \multirow{6}[6]{*}{} & \multicolumn{2}{c}{CPX} & \multirow{6}[6]{*}{} & \multicolumn{3}{c}{OBD} \\
\cmidrule{3-5}\cmidrule{7-8}\cmidrule{10-12}          &       & itr   & min   & gap   &       & min   & gap   &       & itr   & min   & gap \\
\cmidrule{1-5}\cmidrule{7-8}\cmidrule{10-12}    20    & 0.1014  & 9     & 17.5  & \textcolor[rgb]{ 0,  0,  0}{} &       & 399.4  & \textcolor[rgb]{ 0,  0,  0}{} &       & 7     & 22.1  & \textcolor[rgb]{ 0,  0,  0}{} \\
    40    & 0.0953  & 10    & 43.6  & \textcolor[rgb]{ 0,  0,  0}{} &       & T     & N/A   &       & 7     & 163.1  & \textcolor[rgb]{ 0,  0,  0}{} \\
    60    & 0.0863  & 11    & 133.6  & \textcolor[rgb]{ 0,  0,  0}{} &       & T     & N/A   &       & /     & T     & 429.1\% \\
    80    & 0.0857  & 8     & 138.5  & \textcolor[rgb]{ 0,  0,  0}{} &       & T     & N/A   &       & /     & T     & 671.2\% \\
    \bottomrule
    \end{tabular}%
  \label{tab:cp}%
\end{table}%

Furthermore, the comparison between SBD and benchmark algorithms under different numbers of period ($NP$) are shown in Table \ref{tab:ps}. It indicates that the robustness level of solution evidently improves with the increase of $NP$, which, however, makes it more challenging to solve the problem. In all instances, only SBD and OBD can converge to an optimal solution. When $NP=3$, CPX is terminated at the gap of 5.06\% due to the time limit, which leads to a sub-optimal solution. Nevertheless, when $NP=5$, CPX fails to report a feasible solution. In contrast, our SBD demonstrates a remarkable solution capacity, which not only solves the problem successfully, but also saves the solution time by 73.27\%, compared to OBD. Note that when $NP=10$, only SBD can successfully solve the problem, while OBD can only report a low quality solution with a very large gap. Hence, comparing to CPX or OBD, SBD performs faster by orders of magnitude.
\renewcommand\arraystretch{0.7}
\begin{table}[!t]
  \centering
  \setlength{\tabcolsep}{1.4mm}
  \caption{Simple Microgrid: Solution Performance under Different Period Numbers}
    \begin{tabular}{cccccccccccc}
    \toprule
    \multirow{2}[4]{*}{\textit{NP}} & \multirow{2}[4]{*}{OBJ} & \multicolumn{3}{c}{\#SBD\#} & \multirow{5}[6]{*}{} & \multicolumn{2}{c}{CPX} &       & \multicolumn{3}{c}{OBD} \\
\cmidrule{3-5}\cmidrule{7-8}\cmidrule{10-12}          &       & itr   & min   & gap   &       & min`  & gap   &       & itr   & min   & gap \\
\cmidrule{1-5}\cmidrule{7-12}    3     & 0.0055 & 5     & 12    & \textcolor[rgb]{ 0,  0,  0}{} &       & T     & 5.06\% &       & 5     & 37.7  & \textcolor[rgb]{ 0,  0,  0}{} \\
    5     & 0.0953 & 10    & 43.6  & \textcolor[rgb]{ 0,  0,  0}{} &       & T     & N/A   &       & 7     & 163.1 & \textcolor[rgb]{ 0,  0,  0}{} \\
    10    & 0.2295 & 10    & 146   & \textcolor[rgb]{ 0,  0,  0}{} &       & T     & N/A   &       & /     & T     & 378.9\% \\
    \bottomrule
    \end{tabular}%
  \label{tab:ps}%
\end{table}%

\subsection{Complex Microgrid}
\label{Sec:Case2}
The proposed planning method is further applied to optimize the investment portfolio of a practical microgrid with complex structure. Detailed parameters of candidate DERs, including 4 types of RES, 7 types of DFGs, and 4 types of ESS, are listed in Table \ref{tab:der2}. In this case, 4 planning periods are considered over 12 years. The initial load capacity is 1.25MW and the forecasted demand growth is 10\% per year. The budget level limit is specified as $1.3\times \$15,624,525.12$ while the risk tolerance level is set to $0.10$. Table \ref{tab:hmsp} presents the planning result of complex microgrid. We observe that the DERs with lower capital cost (WT2, PV2, and ES4) or lower fuel cost (DE1 and MT1) are installed with a high priority, which ensures the cost-effectiveness of solution. The robustness level of CC-IGD solution is maximized as $\alpha_{\rm L}=0.1052$ by properly sizing the DERs' combination. On one hand, each kilowatt of RES is equipped with 1.23kWh storage device to deal with resource uncertainties. On the other hand, over 2900kW DFGs (30\% of total generation capacity) are installed over 12 years to ensure the power-supply reliability under demand growth uncertainty. Hence, the derived planning scheme gives a full consideration to the complex uncertain environment of microgrid.

\renewcommand\arraystretch{0.7}
\begin{table}[!t]
  \centering
  \caption{Complex Microgrid: Parameters of Candidate DERs}
  \setlength{\tabcolsep}{0.30mm}
    \begin{tabular}{ccccccccc}
    \toprule
    \multicolumn{9}{c}{Distributed Generators} \\
    \midrule
    \multirow{3}[2]{*}{Label} & \multirow{3}[2]{*}{Type} & Rated & Min & Capital & O\&M  & Fuel  & Num &  \\
          &       & Power & Output & Cost  & Cost  & Cost  & Limit &  \\
          &       & (kW)  & (kW)  & (\$/kW) & (\$/kW/yr) & (\$/kWh) &       &  \\
    \midrule
    WT1   & Wind Turbine & 240   & 0     & 900   & 0     & /     & 4     &  \\
    \midrule
    WT2   & Wind Turbine & 360   & 0     & 750   & 0     & /     & 4     &  \\
    \midrule
    PV1   & Photovolatic & 100   & 0     & 1150  & 0     & /     & 4     &  \\
    \midrule
    PV2   & Photovolatic & 180   & 0     & 1000  & 0     & /     & 4     &  \\
    \midrule
    DE1   & Diesel Engine & 120   & 20    & 280   & 131.4 & 0.095 & 2     &  \\
    \midrule
    DE2   & Diesel Engine & 150   & 30    & 250   & 131.4 & 0.099 & 2     &  \\
    \midrule
    DE3   & Diesel Engine & 200   & 40    & 230   & 131.4 & 0.109 & 2     &  \\
    \midrule
    MT1   & Micro Turbine & 160   & 15    & 420   & 236.5 & 0.086 & 2     &  \\
    \midrule
    MT2   & Micro Turbine & 200   & 20    & 380   & 236.5 & 0.091 & 2     &  \\
    \midrule
    FC1   & Fuel Cell & 80    & 0     & 620   & 289.1 & 0.137 & 2     &  \\
    \midrule
    FC2   & Fuel Cell & 100   & 0     & 590   & 289.1 & 0.156 & 2     &  \\
    \midrule
    \multicolumn{9}{c}{Energy Storage Devices} \\
    \midrule
    \multirow{3}[2]{*}{Label} & \multirow{3}[2]{*}{Type} & Rated & Rated & Power & Energy  & O\&M  & Effi- & Num \\
          &       &  Power & Energy & Cost  & Cost  & Cost  & \multirow{2}[1]{*}{ciency} & \multirow{2}[1]{*}{Limit} \\
          &       & (kW)  & (kWh) & (\$/kW) & (\$/kWh) & (\$/kW/yr) &       &  \\
    \midrule
    ES1   & Battery Pack 1 & 100   & 200   & 340   & 220   & 20    & 0.95  & 6 \\
    \midrule
    ES2   & Battery Pack 2 & 150   & 200   & 280   & 220   & 20    & 0.9   & 6 \\
    \midrule
    ES3   & Battery Pack 3 & 200   & 300   & 250   & 170   & 30    & 0.85  & 6 \\
    \midrule
    ES4   & Battery Pack 4 & 300   & 300   & 200   & 170   & 30    & 0.8   & 6 \\
    \bottomrule
    \end{tabular}%
  \label{tab:der2}%
\end{table}%

Also, the computational test of proposed solution approach is performed on a complex microgrid under different numbers of scenarios. The solution time limit is set to 2000 mins, which is moderate for a practical planning problem. As indicated by the test results in Table \ref{tab:hmcp}, this practical system with 10 scenarios is extremely challenging for CPX to compute. This instance, however, can be well addressed by both bilinear BD methods. Moreover, with the number of scenarios growing, SBD demonstrates its superiority over OBD by successfully solving all the instances (including an intractable 80-scenario instance) using at most 1844 mins. On the contrary, OBD fails to close the optimality gap for instances with more than 10 scenarios. Together with our observations made on the computational study using the simple microgrid, we can conclude that the customized SBD has a superior scalable capacity to handle MMEP with multiple planning periods and stochastic scenarios.

\renewcommand\arraystretch{0.7}
\begin{table}[!t]
  \centering
  \caption{Complex Microgrid: Planning Results}
  \setlength{\tabcolsep}{2.2mm}
    \begin{tabular}{ccc}
    \toprule
    Period & DER Type & Planning Scheme \\
    \midrule
    \multirow{5}[6]{*}{I
Year 1-3} & \multirow{2}[2]{*}{RES} & WT1=480kW WT2=1440kW  \\
          &       & PV1=400kW PV2=720kW \\
\cmidrule{2-3}          & \multirow{2}[2]{*}{DFG} & DE1=240kW DE2=300kW \\
          &       &  MT1=320kW MT2=400kW \\
\cmidrule{2-3}          & BES   & ES3=1800kWh ES4=1800kWh \\
    \midrule
    \multirow{3}[6]{*}{II
Year 4-6} & RES   & WT2=1440kW PV2=720kW \\
\cmidrule{2-3}          & DFG   & DE1=240kW DE2=150kW MT1=320kW \\
\cmidrule{2-3}          & BES   & ES3=1800kWh ES4=1800kWh \\
    \midrule
    \multirow{3}[6]{*}{III
Year 7-9} & RES   & WT2=1080kW PV2=180kW \\
\cmidrule{2-3}          & DFG   & DE1=240kW MT1=320kW \\
\cmidrule{2-3}          & BES   & ES4=1200kWh \\
    \midrule
    \multirow{2}[4]{*}{IV
Year 10-12} & RES   & WT2=360kW \\
\cmidrule{2-3}          & DFG   & DE1=120kW MT1=320kW \\
    \bottomrule
    \end{tabular}%
  \label{tab:hmsp}%
\end{table}%

\renewcommand\arraystretch{0.7}
\begin{table}[!t]
  \centering
  \setlength{\tabcolsep}{1.4mm}
  \caption{Complex Microgrid: Solution Performance of Different Methods}
    \begin{tabular}{cccccccccccc}
    \toprule
    \multirow{2}[4]{*}{\textit{NS}} & \multirow{2}[4]{*}{OBJ} & \multicolumn{3}{c}{\#SBD\#} & \multirow{5}[6]{*}{} & \multicolumn{2}{c}{CPX} & \multirow{5}[6]{*}{} & \multicolumn{3}{c}{OBD} \\
\cmidrule{3-5}\cmidrule{7-8}\cmidrule{10-12}          &       & itr   & min   & gap   &       & min   & gap   &       & itr   & min   & gap \\
\cmidrule{1-5}\cmidrule{7-8}\cmidrule{10-12}    10    & 0.1101 & 18    & 115.5 & \textcolor[rgb]{ 0,  0,  0}{} &       & T     & 9.84\% &       & 15    & 613.0  & \textcolor[rgb]{ 0,  0,  0}{} \\
    40    & 0.0916 & 16    & 646.6 &       &       & T     & N/A   &       & /     & T     & 315.2\% \\
    80    & 0.1052 & 14    & 1844.0  & \textcolor[rgb]{ 0,  0,  0}{} &       & T     & N/A   &       & /     & T     & 689.3\% \\
    \bottomrule
    \end{tabular}%
  \label{tab:hmcp}%
\end{table}%

\section{Conclusion}
\indent In this paper, a chance constrained IGD model is proposed to manage the uncertainties in MMEP problem. This model is formulated under an integrated framework of IGDT and chance constrained program, which aims to maximize the robustness level of DER investment meanwhile satisfying a set of operational constraints with a high probability. Furthermore, a strengthened bilinear Benders decomposition algorithm is developed to solve the problem. Two sets of experiments are presented to verify the proposed planning method. We observe that our chance constrained IGD model is competent to address the risks and challenges from both random and non-random uncertainties within a multi-period planning scheme. Also, our strengthened bilinear Benders Decomposition method demonstrates a strong solution capacity and scalability to compute the chance constrained IGD model for practical microgrids. Therefore, it enables us to make informative planning decision in a more efficient way.

%







\ifCLASSOPTIONcaptionsoff
  \newpage
\fi



%




\bibliographystyle{IEEEtran}
\bibliography{test1}

%




\end{document}